%% file: polycubes.tex
\documentclass[10pt]{article}

\usepackage[latin1]{inputenc}
\usepackage{subfigure}
\usepackage[english]{babel}
\usepackage[T1]{fontenc}
\usepackage{float}
\usepackage[margin=0.7in]{geometry}

\usepackage{wrapfig}
\restylefloat{figure}
\usepackage{amsmath}
\usepackage{cite}
\usepackage{amsmath,amssymb}
\usepackage{lipsum}

\usepackage{pgf}
\usepackage{tikz}
\usetikzlibrary{arrows,automata,plothandlers,positioning,calc}

\usepackage[normalem]{ulem}
\usepackage[margin=10cm]{caption}

\newtheorem{theorem}{Theorem}
\newtheorem{cor}[theorem]{Corollary}
\newtheorem{lem}[theorem]{Lemma}
\newtheorem{prop}[theorem]{Proposition}

\def\N{{\mathbb N}}
  
\def\cqfd{\hfill $\Box$ \bigskip}

\newcommand{\cellule}[2]{
			\fill[rounded corners, color=black!50] (#1,#2) +(.02,.02
) -- +(.02,.98) -- +(.98,.98) -- +(.98,.02) -- cycle;
		}
\newcommand{\rcellule}[2]{
			\fill[rounded corners, color=red] (#1,#2) +(.02,.02
) -- +(.02,.98) -- +(.98,.98) -- +(.98,.02) -- cycle;
		}

\newcommand{\gcellule}[2]{
			\fill[rounded corners, color=green] (#1,#2) +(.02,.02
) -- +(.02,.98) -- +(.98,.98) -- +(.98,.02) -- cycle;
		}

\newcommand{\cube}[3]{
			\fill[color = black!20] (#1,#2,#3) +(0,1,0) -- +(0,1,-1)
 -- +(1,1,-1) -- +(1,1,0) -- cycle;
			\fill[color = black!30] (#1,#2,#3) -- +(0,1,0) -- +(1,1,
0) -- +(1,0,0) -- cycle;
			\fill[color = black!50] (#1,#2,#3) +(1,0,0) -- +(1,1,0) 
-- +(1,1,-1) -- +(1,0,-1) -- cycle;
			\draw[color = black!80] (#1,#2,#3) -- +(0,1) -- +(1,1) -
- +(1,0) -- cycle;
			\draw[rounded corners,color = black!80] (#1,#2,#3) + (0
.025,.94,0) -- +(0,1,-1) -- +(1,1,-1) -- +(1,0,-1) -- +(1,0,0);
			\draw[color = black!80] (#1,#2,#3) + (0,1,0) -- +(0,1,-1
) -- +(1,1,-1) -- +(1,0,-1) -- +(1,0,0);
			\draw[color = black!80] (#1,#2,#3) +(1,1,0) -- +(1,1,-1)
;
}

\newcommand{\tcube}[3]{
			\draw[color = black!20] (#1,#2,#3) +(0,1,0) -- +(0,1,-1)
 -- +(1,1,-1) -- +(1,1,0) -- cycle;
			\draw[color = black!30] (#1,#2,#3) -- +(0,1,0) -- +(1,1,
0) -- +(1,0,0) -- cycle;
			\draw[color = black!50] (#1,#2,#3) +(1,0,0) -- +(1,1,0) 
-- +(1,1,-1) -- +(1,0,-1) -- cycle;
			\draw[color = black!80] (#1,#2,#3) -- +(0,1) -- +(1,1) -
- +(1,0) -- cycle;
			\draw[rounded corners,color = black!80] (#1,#2,#3) + (0
.025,.94,0) -- +(0,1,-1) -- +(1,1,-1) -- +(1,0,-1) -- +(1,0,0);
			\draw[color = black!80] (#1,#2,#3) + (0,1,0) -- +(0,1,-1
) -- +(1,1,-1) -- +(1,0,-1) -- +(1,0,0);
			\draw[color = black!80] (#1,#2,#3) +(1,1,0) -- +(1,1,-1)
;
	\draw[color = black!80] (#1,#2,#3-1) -- +(0,1) -- +(1,1) -
- +(1,0) -- cycle;
	\draw[color = black!80] (#1,#2,#3) +(0,0,0) -- +(0,1,0) 
-- +(0,1,-1) -- +(0,0,-1) -- cycle;
}

\newcommand{\bcube}[3]{
			\fill[color = red] (#1,#2,#3-1) -- +(0,1) -- +(1,1) -
- +(1,0) -- cycle;
}

\newcommand{\lcube}[3]{
			\fill[color = green] (#1,#2,#3) +(0,0,0) -- +(0,1,0) 
-- +(0,1,-1) -- +(0,0,-1) -- cycle;
;
}

\def\C{{\mathbb C}}
\def\Z{{\mathbb Z}}
\def\N{{\mathbb N}}

\def\P{{\mathbb P}}
\def\E{{\mathbb E}}

\author{
C. Carr\'e\and
N. Debroux\and
M. Deneufch\^atel\and
J.-Ph. Dubernard\and
C. Hillairet\and
J.-G. Luque\and
O. Mallet
\thanks{\texttt{christophe.carre,matthieu.deneufchatel,jean-philippe.dubernard,
jean-gabriel.luque,olivier.mallet@univ-rouen.fr}}\ \thanks{\texttt{noemie.debroux,conrad.hillairet@insa-rouen.fr}}}

\title{Dirichlet convolution and enumeration of pyramid polycubes
}
\def\href#1{\texttt{#1}}
\def\url#1{\texttt{#1}}
\begin{document}
\maketitle
\begin{abstract} 
{\bf Abstract:}
We investigate the enumeration of two families of polycubes, namely pyramids and espaliers, in connection with the multi-indexed Dirichlet convolution.\\
{\bf R\'esum\'e:} Nous \'etudions l'\'enum\'eration de deux familles de polycubes, les pyramides et les espaliers, en lien avec une version multi-index\'ee de la convolution 
de Dirichlet.
\end{abstract}

\input{intro_m}

\input{definition}

\input{convolution}

\input{countingpyr}

\input{general}
\input{concl}

%
%
%
%

%

\bibliographystyle{alpha}
\input{ref}

\end{document}

%% file: intro_m.tex
\section{Introduction}

%
%
%
In the Cartesian plane $\mathbb{Z}^2$, a polyomino is a
finite connected union of elementary cells (unit squares) without cut
point and defined up to a translation. Even if they have been studied for a long time in
combinatorics, no exact formula is known for counting general polyominoes but many
results have been found concerning certain classes of 
polyominoes, see for instance \cite{Bo96} or \cite{Fer04}. 

\noindent Polyominoes also have a 3-dimensional equivalent: the
 3-dimensional polycubes (or polycubes for short) \cite{Lu71}.
If we consider, now, that an elementary cell is a unit cube, then a polycube 
is a face-connected finite set of elementary cells defined up to a translation 
in $\mathbb{Z}^3$. Like polyominoes, polycubes appear in statistical physics,
 more 
precisely in the phenomenon of percolation (see \cite{BH57} for example). A lot 
of studies have led to count polycubes with respect to their number, $n$ say, of 
cells. The first values were found in 1972 up to $n=6$ 
\cite{Lu71} and the last one (to our knowledge) in 2006, up to $n=18$ \cite{AB06}.

\noindent The notion of polycube can be extended to dimension $d$, with 
$d\geq 3$; $d$-dimensional polycubes (or $d$-polycubes for short) are used in an efficient model of real-time
validation \cite{LG08}, as well as in the representation of finite geometrical 
languages \cite{CDJ09,Jea10}. Although the polycubes are higher dimensional natural analogues of polyominoes, very little is known about their enumeration. 
In particular only few families of polycubes have been studied. In this paper, we propose to investigate two  classes of 
polycubes: pyramids and espaliers. The interest of these two examples lies in their connection with Lambert and Dirichlet 
generating series.\\
The paper is organised as follows. First, in Section \ref{families}, we define pyramids and espaliers in dimension $d+1$ and 
investigate the first properties. In particular, we show that espaliers of height $h$ make it possible to describe a partial order on
partitions of height $h$ which recovers the classical division order for $h=1$. In Section \ref{convolution}, we extend the
 convolution product to multi-indexed 
 families and we give an interpretation in terms of ordinary and Dirichlet generating functions.
   Furthermore, we point out the connection with espaliers. In Section \ref{counting} we apply the properties
    of the convolution product
   to the enumeration of pyramids and  espaliers. In particular, we show that the number $n_v(d)$ of pyramids of volume $v$ in dimension
   $d+1$ is a polynomial in $d$ of degree $\lfloor \log_2(v)\rfloor$. Finally, in Section \ref{general} we explain how to apply our method to other families of
   polycubes.

%% file: definition.tex
\section{Some families of polycubes\label{families}}
%
%


\subsection{Definitions}

We will consider polycubes as discrete objects which are embedded in the three-dimensional discrete lattice $\Z^3$.
Each point of $\Z^3$ will be represented by the triplet of its coordinates and lexicographically ordered. An atomic cell is a cube of volume $1$ which will be identified with the smallest coordinates of its vertices.
So, for our purpose and without loss of generality, we will consider a polycube as a finite and connected (by face) collection $\mathcal P$ of cells such that
its smallest cell is $(0,0,0)$.
The \emph{volume} of a polycube is the number of its atomic cells and its
\emph{height}  is the difference between the 
greatest and the 
smallest indices of its cells  according to the first coordinate.

\noindent 
A very particular  polycube is the \emph{horizontal plateau}: it is a horizontal parallelepiped of height 1. 
To simplify the notations, let us call it a plateau. 
The notion of plateau allows us to define two new families of polycubes. They appear in the study of two particular families of convex-directed polycubes \cite{CohenSolal}. The first family is a subclass of plane partitions (see \cite{CLP98,Bres99,Stan01} for instance).

A \emph{pyramid} polycube (or  \emph{pyramid} for short) is obtained 
by gluing together horizontal plateaus in such a way that
\begin{itemize}
\item $(0,0,0)$ belongs to the first plateau, and each cell with  coordinates
 $(0,b,c)$ belonging to the first plateau is such that
$b,c\geq 0$.
\item If the cell  with coordinates $(a,b,c)$ belongs to the $(a+1)$-th plateau ($a>0$), 
then the cell with  coordinates $(a-1, b, c)$ belongs to the $a$-th plateau.
\end{itemize}

   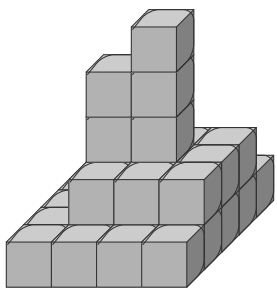
\begin{wrapfigure}{l}{0.25\textwidth}
 
 \begin{minipage}{\marginparwidth}%
 \centering
		\begin{tikzpicture}[scale=0.6]
				
				\begin{scope}[xshift = 0cm]
					\foreach \x in {0,1,2,3}{
					\foreach \z in {0,1,2,3,4}{
						\cube{\x}{0}{\z}
					
}					}

				\foreach \x in {1,2,3}{
					\foreach \z in {1,2,3}{
						\cube{\x}{1}{\z}
					
}					}

				\foreach \x in {1,2}{
					\foreach \z in {2}{
						\cube{\x}{2}{\z}
					
}					}
				\foreach \x in {1,2}{
					\foreach \z in {2}{
						\cube{\x}{3}{\z}
					
}					}
\cube{2}{4}{2}
				\end{scope}
			\end{tikzpicture}
    \end{minipage}	
    \caption{A pyramid}
	\label{pyramid1}
	\end{wrapfigure}
An \emph{espalier} polycube is a special pyramid such that each plateau contains the cell $(a,0,0)$.\\
\subsection{Counting pyramids}
The most natural statistic to count   pyramids and   espaliers is the volume. 
The number of pyramids and espaliers of a given volume are presented below up to volume 12:
\[
1,3,7,16,33,63,117,202,344,566,908,1419
\]
\[
1,3,5,10,14,26,34,57,76,116,150,227
\] 
They correspond respectively to the sequences \href{http://oeis.org/A229914}{A229914} and \href{http://oeis.org/A229915}{A229915} in \cite{Sloane}.
This statistic can be refined
by the height (\emph{i.e.} the number of plateaus) and the volume of each plateau. 
If $\mathbb P$ is a pyramid of height $h$, we will denote by $mv(\mathbb P)=(v_1,\dots,v_h)$ with
$v_1\geq v_2\geq\dots \geq v_h>0$ the sequence of the volumes of its plateaus. 
Let $\lambda=(\lambda_1\geq\cdots\geq\lambda_h)$, with $\lambda_h>0$ be a partition, 
we define $E_\lambda$ as the set of espaliers $\E$ such that
$mv(\E)=\lambda$.
We use the following notations where $t=e$ if the corresponding quantity involves espalier polycubes and $t=p$ if it involves pyramid polycubes:
\begin{itemize}
 \item the number of objects of volume $v$  is denoted by $n_v^{t}$; we denote by $n_{v,h}^{t}$ 
 the number of objects of given height $h$ and volume $v$;
 \item there are $n^t_{\left[ v_1 , \dots , v_h \right]}$ objects such that each plateau has volume 
 $v_i$, $1 \leq i \leq h$, if $v_1 \geq \dots \geq v_h$.
\end{itemize}
By convention, there is no espalier nor pyramid of volume $0$: $n^{t}_0= 0$. 
The number $n^t_{i,j,h,v}$ of considered polycubes (espaliers or pyramids) of volume $v$, height $h$ and such that its largest plateau is $i\times j$ is given by the recurrence:
$
n_{i,j,h,v}^t=
	\sum_{a,b}\alpha^t_{a,b}n_{i+a,j+b,h-1,v-ij}^t $
 with $\alpha^e_{a,b}=1,\alpha^p_{a,b}=(a+1)(b+1)$ and $n^t_{i,j,1,v}=\delta_{ij,v}$.


 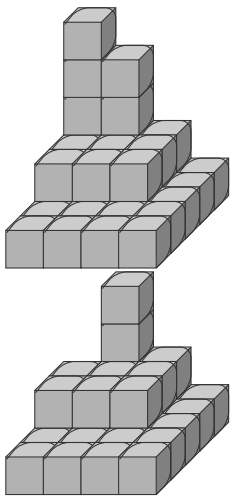
\begin{wrapfigure}{l}{0.30\textwidth}
 
 \centering
		\begin{tikzpicture}[scale=0.5]
				
				
				\begin{scope}[xshift = 0cm]
					\foreach \x in {0,1,2,3}{
					\foreach \z in {0,1,2,3,4}{
						\cube{\x}{0}{\z}
					
}					}

				\foreach \x in {0,1,2}{
					\foreach \z in {0,1,2}{
						\cube{\x}{1}{\z}
					
}					}

				\foreach \x in {0,1}{
					\foreach \z in {0}{
						\cube{\x}{2}{\z}
					
}					}
				\foreach \x in {0,1}{
					\foreach \z in {0}{
						\cube{\x}{3}{\z}
					
}					}
\cube{0}{4}{0}
				\end{scope}
			\end{tikzpicture}
 		\begin{tikzpicture}[scale=0.5]
				
				
				\begin{scope}[xshift = 0cm]
					\foreach \x in {0,1,2,3}{
					\foreach \z in {0,1,2,3,4}{
						\cube{\x}{0}{\z}
					
}					}
\cube{0}{1}{1}
\cube{0}{1}{2}

				\foreach \x in {1,2}{
					\foreach \z in {0,1,2}{
						\cube{\x}{1}{\z}
					
}					}

				\foreach \x in {1}{
					\foreach \z in {0}{
						\cube{\x}{2}{\z}
					
}					}
				\foreach \x in {1}{
					\foreach \z in {0}{
						\cube{\x}{3}{\z}
					
}					}

				\end{scope}
			\end{tikzpicture}
 
    \caption{An espalier and a its associated quasi-espalier}
	\label{espalier1}
	\end{wrapfigure}

The generating function of the number of considered polycubes (espaliers or pyramids) of given height $h$ is denoted by 
${\cal E}^t(x;h) := \sum_{v \geq 1} n_{v,h}^t x^v$; we define also 
${\cal E}^t(x) := \sum_{v \geq 1} n_{v}^t x^v=\sum_h {\cal E}^t(x;h)$. Taking into account the distribution of
 the volume among the levels, one gets ${\cal E}^t (x_1 , \dots , x_h ; h) = 
 \sum_{m_1 \geq \dots \geq m_h} n^t_{\left[ m_1 , \dots , m_h \right]} x_1^{m_1} \dots x_h^{m_h}$. 

We also present results about Dirichlet generating functions which are defined by:
\[
 {\cal E}_{\cal D}^t (s_1 , \dots , s_h ; h) = \sum_{m_1 \geq \dots \geq m_h} \frac{n^t_{[m_1,\dots,m_h]}}
{m_1^{s_1} \dots m_h^{s_h}}.
\]
It is interesting to note that the limit $\lim\limits_{h \to \infty} x^{-h} {\cal E}^e(x;h)$ exists. 
The corresponding series is in fact the generating function of a class of polycubes, which we call \emph{quasi-espaliers}, counted by volume. 
 Quasi-espaliers are espaliers from which all the cells with coordinates $(a,0,0)$ have been removed. Note
 that quasi-espaliers are not considered up to a translation: they are the figures obtained when we remove the column $(a,0,0)$
 from an espalier based at $(0,0,0)$ (see Figure \ref{espalier1} for an example). For instance $\{(0,1,0)\}$ and $\{(0,0,1)\}$ are different quasi-espaliers. 
 The first values are
$
2, 4, 7, 12, 18, 29, 42, 61, 87, 122, 167, 229.
$

If we extend the notion of quasi-espalier to pyramids, we can also extend the previous result to pyramids. A \emph{ quasi-pyramid} of is obtained from a pyramid of height $h$ by choosing a cell $(h,b,c)$ in the pyramid and deleting all the cells $(a,b,c)$ with $1\leq a\leq h$. 
Let $Q^p(x)$ be their generating function with respect to volume. Then $\lim_{h\rightarrow \infty}x^{-h}\mathcal E^p(x,h)$ exists and $\lim_{h\rightarrow \infty}x^{-h}\mathcal E^p(x,h)=Q^p(x)+{x\over 1-x}.$ 
 
\subsection{The projection order\label{order}}
    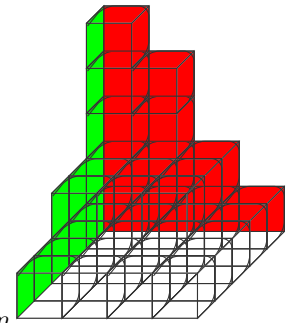
\begin{wrapfigure}{r}{0.30\textwidth}
 
 \centering
		\begin{tikzpicture}[scale=0.6]
				
				\begin{scope}[xshift = 0cm]
				\lcube{0}{0}{0}\lcube{0}{0}{1}\lcube{0}{0}{2}\lcube{0}{0}{3}\lcube{0}{0}{4}
\lcube{0}{1}{0}\lcube{0}{1}{1}\lcube{0}{1}{2}
\lcube{0}{2}{0}
\lcube{0}{3}{0}
\lcube{0}{4}{0}
				\foreach \x in {0,1,2,3}{
					\bcube{\x}{0}{0}
					\foreach \z in {0,1,2,3,4}{
						\tcube{\x}{0}{\z}
					
}					}

				\foreach \x in {0,1,2}{
					\bcube{\x}{1}{0}
					\foreach \z in {0,1,2}{
						\tcube{\x}{1}{\z}
					
}					}
				\foreach \x in {0,1}{
				\bcube{\x}{2}{0}
					\foreach \z in {0}{
						\tcube{\x}{2}{\z}
					
}					}
				\foreach \x in {0,1}{
				\bcube{\x}{3}{0}
					\foreach \z in {0}{
						\tcube{\x}{3}{\z}
					
}					}
\bcube{0}{4}{0}\tcube{0}{4}{0}

				\end{scope}
			\end{tikzpicture}
    \caption{One-to-one correspondence between espaliers and pairs of partitions}
	\label{bijection1}
	\end{wrapfigure}
We describe here an order on integer partitions of same height $h$. 
Note first, that there is a one-to-one correspondence between
espaliers of height $h$ and pairs of partitions of the same height $h$.
 The partitions are obtained by projecting an espalier $\mathbb E$ of height $h$ as
$\mathbb E_x:=\{(a,0,c):(a,b,c)\in\mathbb E\}$ and  $\mathbb E_y:=\{(a,b,0):(a,b,c)\in\mathbb E\}.$
These two sets are obviously two Ferrers diagrams which represent two partitions $\lambda_x(\mathbb E)$ and 
$\lambda_y(\mathbb E)$ of height $h$ (see Fig. \ref{bijection1} for an example) and the bijection is straightforward from the construction.\\
We define the relation $\preceq$ on the set of partitions of height $h$ by $\lambda\preceq\mu$ if and only if there exists 
$\E\in E_\mu$ such that $\lambda_x(\E)=\lambda$.
\begin{prop}
The relation $\preceq$ is a partial order which generalizes the division order on integers in the following sense:
$(d,\dots,d)\preceq (n,\dots,n)\mbox{ if and only if } d|n.$
\end{prop}
Although it is not the purpose of this article, we note that the M\"obius function $\mu^h$ of this order 
seems to have interesting properties. For instance for $h=2$, it satisfies the following equalities:
\begin{eqnarray*}
\mu^2((1,1),(n,n))=\mu^2((1,1),(n,1))=\mu(n),\, \sum_{k=1}^n\mu((1,1),(n,k))=n\mu(n),\\
\mu^2((1,1),(p,m))=-1 \mbox{ for }p\mbox{ prime and }p>m,\end{eqnarray*}
\[
\mu^2((1,1),(pq,n))=2^{k+1}-1\mbox{ when }p, q\mbox{ are prime and }n\mbox{ is the product of }k\mbox{ prime distinct integers,}\]
\[\mu^2((1,1),(n^2,m^2))=2\mbox{ if }n\geq m^2\mbox{ and }0 \mbox{ otherwise}.\] 
\subsection{Higher dimensional polycubes}
The $d$-polycubes are a natural extension of the notion of polyomino to dimension $d$, with $d\geq 3$ 
(see \cite{CohenSolal,ChamparnaudDJ13} for instance).
An atomic $d$-cell is a cube of volume $1$ identified with the smallest coordinates of its vertices in $\Z^d$.
 A $d$-polycube
is then a $d$-face-connected finite set of elementary cells, defined up to  translation.
The volume of a $d$-polycube is the number of its elementary cells.
\noindent A \emph{$d$-parallelepiped} is a $d$-polycube $P$ such that 
for some $(n_1,\cdots , n_d)\in \mathbb{N}^d$, any cell with
 coordinates $(\alpha_1,\cdots , \alpha_d)$ satisfying $0\leq\alpha_k\leq n_k$, $1\leq k\leq d$, belongs to $P$.
 Then, a $d$-plateau is a $d$-parallelepiped of height $1$ (that is composed of cells ot the form $(a,n_1,\cdots , n_d)\in \mathbb{N}^d$ for a fixed $a$).\\
A $(d+1)$-\emph{pyramid} is a $(d+1)$-polycube  obtained 
by gluing together  $(d+1)$-plateaus in such a way that
\begin{itemize}
\item the cell $(0,0,\dots,0)$ belongs to the first plateau and each cell with coordinates $(0,n_1,\dots,n_d)$
 belonging to the first plateau is such that
$n_1,\dots,n_d\geq 0$.
\item if the cell with coordinates $(n_0,n_1,\dots,n_d)$ belongs to the $(n_0+1)$-th plateau ($n_0>0$), 
then the cell with coordinates $(n_0-1, n_1,\dots,n_d)$) belongs to the $n_0$-th plateau.
\end{itemize}
A $(d+1)$-\emph{espalier} is a $(d+1)$-pyramid such that each plateau contains the cell $(n_0,0,\dots,0)$.
As for $3$-polycubes, we define : \begin{itemize}
 \item the number of objects of volume $v$, denoted by $n_v^{t}(d)$; we denote by $n_{v,h}^{t}(d)$ 
 the number of objects of given height $h$ and volume $v$;
 \item  the number of  objects such that each plateau has volume 
 $v_i$, $1 \leq i \leq h$, if $v_1 \geq \dots \geq v_h$, denoted by $n^t_{\left[ v_1 , \dots , v_h \right]}(d)$.
\end{itemize}
By convention, there is no espalier nor pyramid of volume $0$ and one object in dimension $0+1$ :
$n^{t}_0= 0$ and $n_v^{t}(0) = 1$. 
The generating function of the number of espaliers of given height $h$ is 
denoted by ${\cal E}^t(x;h,d) = \sum_{v \geq 1} n_v^t(d) x^v$; we define also 
${\cal E}^t(x;d) := \sum_{v \geq 1} n_{v}^t x^v=\sum_h {\cal E}^t(x;h,d)$. 
Taking into account the distribution of the volume between the levels, one sets $${\cal E}^t (x_1 , \dots , x_h ;h, d) = 
\sum_{m_1 \geq \dots \geq m_h} n^t_{\left[ m_1 , \dots , m_h \right]}(d) x_1^{m_1} \dots x_h^{m_h}.$$
The Dirichlet generating functions will be denoted by:
$ \quad {\cal E}_{\cal D}^t (s_1 , \dots , s_h ; h,d) := \sum_{m_1 \geq \dots \geq m_h} \frac{n^t_{v}(d)}{m_1^{s_1} \dots m_h^{s_h}}.
$

%% file: convolution.tex
\section{Multivariate versions of the Lambert transform\label{convolution}}
\subsection{Convolution and multivariate series}
We consider a natural multidimensional generalization of the Dirichlet convolution.
For each $h \in \N$, we consider the set $\mathcal M_h:=\{(a_{n_1 , \dots , n_h})_{n_1 , \dots , n_h\geq 1}:
a_{n_1 , \dots , n_h}\in\C\}$.
 Let $A = \left( a_{n_1 , \dots , n_h} \right)_{n_1 , \dots , n_h},
  B = \left( b_{n_1 , \dots , n_h} \right)_{n_1 , \dots , n_h}$; we denote by 
  $C = A \star B = \left( c_{n_1 , \dots , n_h} \right)_{n_1 \geq \dots \geq n_h}\in\mathcal M_h$ 
 the \emph{multivariate convolution} of $A$ and $B$  defined by 
$
 c_{n_1 , \dots , n_h} = \sum_{n_i = m_i p_i,\, i=1,\dots h} a_{m_1 , \dots , m_h} b_{p_1 , \dots , p_h}.
$
\begin{prop}
For any $h\in\N$, the product $\star$ is distributive and $\mathbf 1=\left(\delta_{1,n_1}\cdots \delta_{1,n_h}\right)_{n_1 , \dots , n_h}$ is its identity.
Hence, this endows $\mathcal M_h$ with a structure of commutative algebra.
\end{prop}
Denoting by $S_{A}(x_1,\dots,x_h):=\sum_{n_1 , \dots , n_h\geq1}{a_{n_1,\dots,n_h}}x^{n_1}\cdots
x^{n_h}$ the ordinary generating function of $A$ and
 $S_{A}^{\cal D}(s_1,\dots,s_h):=\sum_{n_1 , \dots , n_h\geq1}{a_{n_1,\dots,n_h}\over n_1^{s_1}\cdots n_h^{s_h}}$ 
 its Dirichlet generating function, we observe the following fact:
 \begin{prop}
Let $A, B\in\mathcal M_h$, the three following assertions are equivalent:
\begin{enumerate}
\item $C=A\star B$;
 \item $\displaystyle S_{C} ( x_1 , \dots , x_h) = \sum_{n_1 , \dots , n_h} a_{n_1 , \dots , n_h} S_B (x_1^{n_1} , \dots , x_h^{n_h})$;
\item $
 \displaystyle S^{{\cal D}}_{C} ( s_1 , \dots , s_h) = S^{{\cal D}}_{A} ( s_1 , \dots , s_h) \, S^{{\cal D}}_{B} ( s_1 , \dots , s_h).$
\end{enumerate}
\end{prop}
Indeed, the map $A\rightarrow S_A$ allows us to endow the ideal $x_1\dots x_h\C[x_1,\dots,x_h]$ with a structure of 
commutative algebra $(x_1\dots x_h\C[x_1,\dots,x_h],+,\star)$ isomorphic to $\mathcal M_h$. With these notations, we have $S_{A\star B}=S_A\star S_B$.
The formal substitution $\mathcal D:x_i^n\rightarrow \frac1{n^{s_i}}$ is an isomorphism
from  $(x_1\dots x_h\C[x_1,\dots,x_h],+,\star)$ to the algebra of multivariate Dirichlet formal series in 
the variables $\{s_1,\dots,s_h\}$. Note that the isomorphism above can be realized through 
  an iteration of Mellin transforms (see \emph{e.g.} \cite[Appendix B, Section B.7]{Fla}):
\[
\mathcal D[f]=\frac1{\Gamma(s_1)\dots\Gamma(s_h)}\int_0^\infty\dots\int_0^\infty f\left(e^{-x_1},\dots,e^{-x_h}\right)x_1^{s_1-1}\cdots x_h^{s_h-1}dx_1\dots dx_h
\]
where $\Gamma(s)=\int_0^\infty e^{-x}x^{s-1}dx$ is the  Euler Gamma function. 

Let $\mathcal T_h:=\left\{\left( a_{n_1 , \dots , n_h} \right)_{n_1, \dots,n_h \geq 1}\in\mathcal M_h:
a_{n_1,\dots,n_h}\neq 0\Rightarrow n_1\geq\dots\geq n_h\right\}$. Since $\mathcal T_h$ is stable under the convolution, it is a subalgebra of $\mathcal M_h$.
For simplicity, we will denote by $\left( a_{n_1 , \dots , n_h} \right)_{n_1\geq \dots\geq n_h \geq 1}$ the elements of $\mathcal T_h$.

\subsection{Multivariate Lambert transform}
Let $\triangle = (1)_{m_1 \geq \dots \geq m_h \geq 1}$.
 We call \emph{multivariate Lambert transform} of $A = \left( a_{n_1 , \dots , n_h} \right)_{n_1 \geq \dots \geq n_h\geq1}$ 
 the convolution of $A$ with $\triangle$:
$ {\rm T}_L(A) := A \star \triangle$.
Remark that the (ordinary) generating function of $\triangle$ is
 $\mathcal S_{\triangle}=\frac{x_1 \dots x_h}{(1 - x_1)(1 - x_1 x_2) \dots (1-x_1 \dots x_h)}$ and its Dirichlet
  generating function is
 $\mathcal S_{\triangle}^{\mathcal D}= \mathcal Z(s_1,\dots,s_h)$ where  ${\cal Z}$ 
 denotes the \emph{large multizeta} function
 ${\cal Z}(s_1,\dots,s_h):=\sum_{n_1 \geq \dots \geq n_h\geq1}{1\over n_1^{s_1}\cdots n_h^{s_h}}$
  \cite{cresson,Costermans:2009:NAM}.
The generating function of ${\rm T}_{L}(A)$ is
\[
 S_{{\rm T}_L(A)} = S_\triangle\star S_A=\sum_{n_1 \geq \dots \geq n_h \geq 1} a_{n_1 , \dots , n_h} \frac{x_1^{n_1} \dots x_h^{n_h}}{(1-x_1^{n_1})(1-x_1^{n_1}x_2^{n_2}) \dots (1 - x_1^{n_1} \dots x_h^{n_h})}
\]
and its Dirichlet generating function is given by
\[
 S^{{\cal D}}_{{\rm T}_L(A)}(s_1,\dots,s_h) = {\cal Z}(s_1,\dots,s_h) S^{{\cal D}}_{A} ( s_1 , \dots , s_h).
\]
The multivariate Lambert transform is related to the order defined in Section \ref{order} by the following formula:
\begin{prop}
Setting $\left(\hat{a}_{n_1 , \dots , n_h}\right)_{n_1 , \dots , n_h}:={\rm T}_L(A)$, if $A=\left({a}_{n_1 , \dots , n_h}\right)_{n_1 , \dots , n_h}$, we obtain:
\[
 \hat{a}_{n_1 , \dots , n_h} = \sum_{(m_1 , \dots , m_h) \preceq (n_1 , \dots , n_h)} a_{m_1 , \dots , m_h}.
\]
\end{prop}
As a consequence, the Dirichlet generating function of $\mu^h((1,\dots,1),\lambda)$ is the inverse of $\cal Z$:
\begin{cor}\label{Mob}
\[
 {\cal Z} (s_1 , \dots , s_h)^{-1} = \sum_{(1 , \dots , 1) \preceq (n_1 , \dots , n_h)} \frac{\mu^h((1 , \dots , 1),(n_1 , \dots , n_h))}{n_1^{s_1} \dots n_h^{s_h}}.
\]
\end{cor}
Note that when $h=1$, ${\cal Z}(s)=\zeta(s)$ is the Riemann zeta function and Corollary \ref{Mob} is the classical identity
$\zeta(s)^{-1}=\sum_{n>0}{\mu(n)\over n^s}$.
\subsection{Another transform}
Let $\blacktriangle = \Big( (n_1 - n_2 + 1) \dots (n_{h-1} - n_{h} +1 ) \Big)_{n_1 \geq \dots \geq n_h \geq 1}$. 
Using $\blacktriangle$ and the convolution product defined above, we construct a new transformation:
$
 {\rm T}_\blacktriangle(A) = \blacktriangle \star A = (a^\blacktriangle_{n_1 , \dots , n_h})_{n_1 \geq \dots \geq n_h \geq 1}.
$
\begin{lem}
\label{lemgenfunc}
   The generating function of $\blacktriangle$ is given by the following formula:
\begin{equation}
\begin{split}
S_\blacktriangle:=\displaystyle \sum_{n_1 \geq \dots \geq n_h \geq 1} (n_1 - n_2 + 1) \dots & ( n_{h-1} - n_{h} + 1 ) x_1^{n_1} \dots x_{h}^{n_h} = \\ 
&  \displaystyle \frac{x_1 \dots x_h}{(1 - x_1)^2 (1 - x_1 x_2)^2 \dots ( 1 - x_1 \dots x_{h-1})^2 ( 1 - x_1 \dots x_h)}.
\end{split}
\end{equation}
\end{lem}
A simple induction on the number of variables proves the result.

Lemma \ref{lemgenfunc} implies that the generating function of ${\rm T}_\blacktriangle(A)$ is 
\[
  S_\blacktriangle\star S_A=\sum_{n_1 \geq \dots \geq n_h \geq 1} a_{n_1 , \dots , n_h} 
 \frac{x_1^{n_1} \dots x_h^{n_h}}{(1-x_1^{n_1})^2 (1-x_1^{n_1}x_2^{n_2})^2 \dots 
 (1-x_1^{n_1} \dots x_{h-1}^{n_{h-1}})^2 (1-x_1^{n_1} \dots x_h^{n_h}) }
\]
and its Dirichlet generating function is given by
$
 S^{\cal D}_{{\rm T}_\blacktriangle(A)} = {\cal Z}^\blacktriangle(s_1 , \dots , s_h)
  \, S^{\cal D}_{A} (s_1 , \dots , s_h) 
$
where
$$
 {\cal Z}^\blacktriangle(s_1 , \dots , s_h) = 
 \sum_{n_1 \geq \dots \geq n_h \geq 1} \frac{(n_1 - n_2 + 1) \dots (n_{h-1} - n_{h}+1)}{n_1^{s_1} \dots n_h^{s_h}}.
$$
As a consequence, we have 
$$
 {a}^\blacktriangle_{n_1 , \dots , n_h} = \sum_{(m_1 , \dots , m_h) \preceq (n_1 , \dots , n_h)} \alpha_{m_1 , \dots , m_h}^{n_1 , \dots , n_h} a_{m_1 , \dots , m_h}.
$$ 
where $$\alpha_{m_1 , \dots , m_h}^{n_1 , \dots , n_h}=\left({n_1\over m_1}-
{n_2\over m_2}+1\right)\cdots\left({n_{h-1}\over m_{h-1}}-
{n_h\over m_h}+1\right)$$ are non-negative integers.

%% file: countingpyr.tex
\section{Application to the enumeration of pyramids\label{counting}}
\subsection{Counting pyramids and espaliers by volume}

 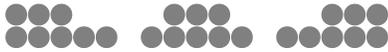
\begin{wrapfigure}{l}{0.6\textwidth}
 
 \centering
		\begin{tikzpicture}[scale=0.3]
\foreach \x in {0,1,2,3,4} {\cellule {\x}{0}}
\foreach \x in {0,1,2,3,4} {\cellule {\x+6}{0}}
\foreach \x in {0,1,2,3,4} {\cellule {\x+12}{0}}
\foreach \x in {0,1,2} {\cellule {\x}{1}}
\foreach \x in {0,1,2} {\cellule {\x+7}{1}}
\foreach \x in {0,1,2} {\cellule {\x+14}{1}}
			\end{tikzpicture}
    \caption{The $3$ $(1+1)$-pyramids of type $(5,3)$}
	\label{(1+1)-pyramids}
	\end{wrapfigure}

The one-to-one correspondence described in Section \ref{order} allows us to construct
an espalier $\mathbb E_{\lambda,\lambda'}\in E_{(\lambda_1\lambda'_1,\dots,\lambda_h\lambda'_h)}$ for each couple of 
partitions $\lambda=(\lambda_1,\dots,\lambda_h)$ and $\lambda'=(\lambda'_1,\dots,\lambda'_h)$. Hence, we deduce
\begin{lem}
\[
\triangle^{\star 2}=\left(n^e_{[m_1,\dots,m_h]}\right)_{m_1\geq\dots\geq m_h\geq 1}.
\]
\end{lem}

   \begin{wrapfigure}{r}{0.25\textwidth}
 
 \begin{minipage}{\marginparwidth}%
 \centering
		\begin{tikzpicture}[scale=0.6]
				
	\lcube{0}{0}{0}\lcube{0}{0}{1}\lcube{0}{0}{2}\lcube{0}{0}{3}\lcube{0}{0}{4}
\lcube{0}{1}{1}\lcube{0}{1}{2}\lcube{0}{1}{3}
\lcube{0}{2}{2}
\lcube{0}{3}{2}
\lcube{0}{4}{2}

	\bcube{0}{0}{0}\bcube{1}{0}{0}\bcube{2}{0}{0}\bcube{3}{0}{0}
\bcube{1}{1}{0}\bcube{2}{1}{0}\bcube{3}{1}{0}
\bcube{1}{2}{0}\bcube{2}{2}{0}
\bcube{1}{3}{0}\bcube{2}{3}{0}
\bcube{2}{4}{0}
				
				\begin{scope}[xshift = 0cm]
					\foreach \x in {0,1,2,3}{
					\foreach \z in {0,1,2,3,4}{
						\tcube{\x}{0}{\z}
					
}					}

				\foreach \x in {1,2,3}{
					\foreach \z in {1,2,3}{
						\tcube{\x}{1}{\z}
					
}					}

				\foreach \x in {1,2}{
					\foreach \z in {2}{
						\tcube{\x}{2}{\z}
					
}					}
				\foreach \x in {1,2}{
					\foreach \z in {2}{
						\tcube{\x}{3}{\z}
					
}					}
\tcube{2}{4}{2}
				\end{scope}
			\end{tikzpicture}
    \end{minipage}	
    \caption{A $(2+1)$-pyramid and its two $(1+1)$-associated pyramids}
	\label{(1+1)-pyramid1}
	\end{wrapfigure}
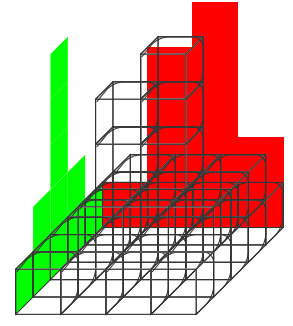
A similar method can be used to compute the number of pyramids. We consider
pyramids in dimension $1+1$. These objects are obtained from partitions (which are espaliers in dimension $1+1$)
by shifting each $(1+1)$-plateau. A $(1+1)$-pyramid $\P$ will be of type $\lambda=(\lambda_1,\dots,\lambda_h)$ if
the first plateau is of size 
$\lambda_1$, $\dots$, the $h$-th plateau is of size $\lambda_h$; we will write $mv(\P)=\lambda$ as with pyramids in dimension $2+1$. 
The number of $(1+1)$-pyramids of type $\lambda$ 
 equals $({\lambda_1}-{\lambda_2}+1)\cdots ({\lambda_{h-1}}-{\lambda_h}+1)$ 
(see Fig. \ref{(1+1)-pyramids} for an example). The pyramids of height $h$ are in one-to-one correspondence with the
pairs of $(1+1)$-pyramids of height $h$.  Furthermore the pyramid $\P$ corresponding to a pair $(\P',\P'')$ 
is such that $mv(\P)=(\lambda'_1\lambda''_1,\dots,\lambda'_h\lambda''_h)$ if $mv(\P')=\lambda'$ and $mv(\P'')=\lambda''$.
We deduce
\begin{lem}
\[
\blacktriangle^{\star 2}=\left(n^p_{[m_1,\dots,m_h]}\right)_{m_1\geq\dots\geq m_h\geq 1}.
\]
\end{lem}
Hence, we are now able to compute the generating function.
\begin{theorem}
The Dirichlet generating functions of espaliers and pyramids are respectively
\begin{equation}\label{EPDir}
{\mathcal E}_{\mathcal D}^e(s_1,\dots,s_h;h)={\cal Z}(s_1,\dots,s_h)^2,\,
{\mathcal E}_{\mathcal D}^p(s_1,\dots,s_h;h)={\cal Z}^\blacktriangle(s_1,\dots,s_h)^2.  
\end{equation}
The ordinary generating functions of the espaliers and the pyramids are respectively
\begin{eqnarray}
{\mathcal E}^e(x_1,\dots,x_h;h)=S_\triangle^{\star2}=\sum_{n_1\geq\dots\geq n_h\geq 1}{x_1^{n_1}\cdots x_h^{n_h}\over
(1-x_1^{n_1})(1-x_1^{n_1}x_2^{n_2})\cdots(1-x_1^{n_1}\cdots x_{h}^{n_h})},\,\label{EOrd}\\
{\mathcal E}^p(x_1,\dots,x_h;h)=S_\blacktriangle^{\star2}=\sum_{n_1\geq\dots\geq n_h\geq 1}{(n_1-n_2+1)\cdots (n_{h-1}-n_{h}+1)x_1^{n_1}\cdots x_h^{n_h}\over
(1-x_1^{n_1})^2(1-x_1^{n_1}x_2^{n_2})^2\cdots(1-x_1^{n_1}\cdots x_{h-1}^{n_{h-1}})^2(1-x_1^{n_1}\cdots x_{h}^{n_h})}.\label{POrd}  
\end{eqnarray}
\end{theorem}
\subsection{Higher dimensions}
In this section, we investigate the coefficients $n_v^t(d)$ ($t=e,p$).
 The iteration of the transformations presented above gives the generating functions of the same objects
 in higher dimension. 
 Hence, the generating function of the number of espaliers (resp. pyramids) whose plateaus have volume $v_1 , \dots , v_h$ 
 in dimension $d+1$ is given by $\triangle^{\star d}$ (resp. $\blacktriangle^{\star d}$).\\
 Therefore the Dirichlet generating function of $n_{[m_1,\dots,m_n]}^t(d)$  is
 \begin{equation}
\mathcal E^t_{\mathcal D}(s_1,\dots,s_h)={\cal Z}^t(s_1,\dots,s_h)^d,
 \end{equation} 
 with ${\cal Z}^e={\cal Z}$ and ${\cal Z}^p={\cal Z}^\blacktriangle$.
We consider now $d$ as a formal parameter. We have:
\[
{\partial^k\over \partial d^k}\mathcal E^t_{\mathcal D}(s_1,\dots,s_h)=
{\cal Z}^t(s_1,\dots,s_h)^d\log\left({\cal Z}^t(s_1,\dots,s_h)\right)^k.
\]
We observe that ${\cal Z}^t(s_1,\dots,s_h)=1+\widetilde{\cal Z}^t(s_1,\dots,s_h)$ where
$\widetilde{\cal Z}^t(s_1,\dots,s_h)=\sum_{n_1\geq 2\atop n_1\geq\dots\geq n_h\geq 1}{(*)\over n_1^{s_1}\dots n_h^{s_h}}$
; here $(*)$ denotes coefficients depending on the value of $t$.\\ It follows that
\begin{equation} \label{eqdif}
{\partial^k\over \partial d^k}\mathcal E^t_{\mathcal D}(s_1,\dots,s_h)=
\sum_{n_1\geq 2^k\atop n_1\geq\dots\geq n_h\geq 1}{(\Box)\over n_1^{s_1}\dots n_h^{s_h}}.
\end{equation}
$(\Box)$ denotes also coefficients. We deduce from (\ref{eqdif}) that $n^t_{[v_1,\dots,v_h]}(d)$, is a polynomial in
$d$ whose degree is at most   $\log_2(v_1)$. Furthermore 
\begin{equation}\label{degn}\deg(n^t_{[2^k,1,\dots,1]}(d))=k.
\end{equation}
 Hence,
\begin{theorem}
The number $n_v^e(d)$ of $(d+1)$-espaliers of volume $v$  and the number 
$n_v^p(d)$ of $(d+1)$-pyramids of volume $v$ are both polynomials in $d$ of degree
$\lfloor \log_2(v)\rfloor.$
\end{theorem}
{\bf Proof}
Since $n_v^t(d)=\sum_{h\geq 1}\sum_{v_1\geq\dots\geq v_h\geq 1}n^t_{[v_1,\dots,v_h]}(d)$, Eq. (\ref{eqdif}) 
implies that  $\deg(n_v^t(d))\leq\lfloor \log_2(v)\rfloor.$
From (\ref{degn}), the inverse inequality holds if 
 $\displaystyle E_{(2^{\lfloor\log_2v\rfloor},\underbrace{1,\dots,1)}_{(v-2^{\lfloor\log_2v\rfloor})\times}}\neq\emptyset.$ 
This is obviously the case since it suffices to consider an espalier such that the volume of the first plateau is $2^{\lfloor\log_2v\rfloor}$ and with
$v-2^{\lfloor\log_2v\rfloor}$ other plateaus consisting of one cell.  
\cqfd

%% file: general.tex
\section{General construction\label{general}}
\begin{wrapfigure}{l}{0.4\textwidth}
 
 \centering
		\begin{tikzpicture}[scale=0.2]
\foreach \x in {0,1,2,3,4} {\gcellule {\x}{9}}
\foreach \x in {7,8} {\gcellule {\x}{9}}

\foreach \y in {0,1,2} {\rcellule {-1}{\y}}
\foreach \y in {5,8} {\rcellule {-1}{\y}}

\foreach \x in {0,1,2,3,4,7,8} {\foreach \y in{0,1,2,5,8}{\cellule {\x}{\y}}}

			\end{tikzpicture}
    \caption{Construction of a stratum}
	\label{Level}
	\end{wrapfigure}
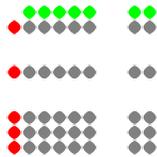

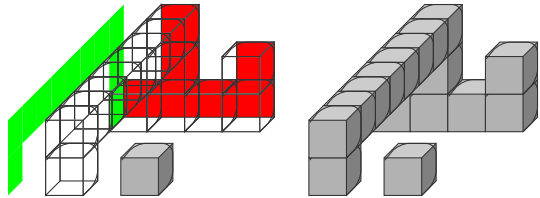
\begin{wrapfigure}{hbtr}{1\textwidth}
 \centering
				



		\begin{tikzpicture}[scale=0.5]
				
	\lcube{0}{0}{0}\lcube{0}{1}{0}\lcube{0}{2}{0}\lcube{0}{2}{1}\lcube{0}{2}{2}
\lcube{0}{2}{3}\lcube{0}{2}{4}\lcube{0}{2}{5}\lcube{0}{2}{6}\lcube{0}{2}{7}\lcube{0}{1}{7}

	\bcube{0}{0}{0}\bcube{1}{0}{0}\bcube{1}{1}{0}\bcube{1}{2}{0}
\bcube{2}{0}{0}\bcube{3}{0}{0}\bcube{3}{1}{0}
				
\tcube{0}{0}{0}\tcube{1}{0}{0}\tcube{2}{0}{0}\tcube{3}{0}{0}
\tcube{1}{1}{0}\tcube{3}{1}{0}
\tcube{1}{2}{0}
\tcube{1}{2}{1}
\tcube{1}{2}{2}
\tcube{1}{2}{3}\tcube{1}{2}{4}\tcube{1}{2}{5}\tcube{1}{2}{6}\tcube{1}{2}{7}
\tcube{1}{1}{7}\cube{3}{1}{7}

\cube{7}{0}{0}\cube{8}{0}{0}\cube{9}{0}{0}\cube{10}{0}{0}
\cube{8}{1}{0}\cube{10}{1}{0}\cube{8}{1}{7}\cube{10}{1}{7}
\cube{8}{2}{0}
\cube{8}{2}{1}
\cube{8}{2}{2}
\cube{8}{2}{3}\cube{8}{2}{4}\cube{8}{2}{5}\cube{8}{2}{6}\cube{8}{2}{7}

			\end{tikzpicture}
    \caption{A counterexample}
	\label{Counterexample}
	\end{wrapfigure}

A $2$-dimensional object is a finite set of cells $(x,y)$ such that $x\in\N$ and $y\in \Z$.
If $\cal O$ is a  $2$-dimensional object, its $i$th stratum will be the set $L_i({\cal O}):=
\{(i,y)\in{\cal O}\}$.
From a pair of $2$-dimensional objects $({\cal O}_1,{\cal O}_2)$, we construct a 
$3$-dimensional object that is a set of cells $(x,y,z)$ with $x\in\N$ and $y,z\in\Z$:
${\cal G}({\cal O}_1,{\cal O}_2):=\{(x,y,z):(x,y)\in{\cal O}_1, (x,z)\in{\cal O}_2)$. The 
$i$-th stratum of a $3$-dimensional object will be
defined by  $L_i({\cal O}):=
\{(i,y,z)\in{\cal O}\}$. The height of an object is $h({\cal O})=\max\{i:L_i({\cal O})\neq \emptyset\}$. 
The multivolume $mv({\cal O})$ of a $2$- or $3$-dimensional object ${\cal O}$
is the sequence of the cardinals of its strata.\\
Note that if $mv({\cal O}_1)=[v_1,\cdots,v_{h}]$ and 
$mv({\cal O}_2)=[v'_1,\cdots,v'_{h'}]$ then $$mv({\cal G}({\cal O}_1,{\cal O}_2))= [v_1v'_1,\dots,
v_{\min(h,h')}v'_{\min(h,h')}].$$ This property is obtained easily by examining the construction of
each stratum (see Fig. \ref{Level} for an example). An object $\cal O$ is said to be plain
if for any $1\leq i\leq h({\cal O})$, the stratum $L_i({\cal O})$ is not empty. 
Let $\mathcal A=\bigcup_{u}\mathcal A_u$ and 
$\mathcal B=\bigcup_{v}\mathcal B_v$ be two families of $2$-dimensional plain objects of height $h$ graded by  
multivolume and such that
$\mathcal A_u$ and $\mathcal B_v$ are finite for any multivolume $u$ and $v$. Set also 
$A=\left(a_{v_1,\dots,v_h}\right)_{v_1,\dots,v_h}$ and  $B=\left(b_{v_1,\dots,v_h}\right)_{v_1,\dots,v_h}$
with $a_{v_1,\dots,v_h}=\#A_{[v_1,\dots,v_h]}$ and $b_{[v_1,\dots,v_h]}=\#B_{[v_1,\dots,v_h]}$.
The set 
${\mathcal G}({\mathcal A},{\mathcal B})={\cal G}({\cal O}_1,{\cal O}_2):
\{{\cal O}_1\in {\mathcal A}_u, {\cal O}_2\in {\mathcal B}_u\}$ contains only plain objects of size $h$. 
Furthermore it is graded
by multivolume ${\mathcal G}({\mathcal A},{\mathcal B})=\bigcup_v{\mathcal G}_v$ and the sequence 
$G=\left(\#G_{[v_1,\dots,v_h]}\right)_{v_1,\dots,v_h\geq 1}$ is given by
$
G=A\star B.
$
Now we want to apply our method in order to construct families of polycubes from families of $2$-dimensional objects.
Note that if $\P_1$ and $\P_2$ are two polyominoes then $G(\P_1,\P_2)$ is not necessarily a polycube. For instance see a
counterexample in Fig \ref{Counterexample}; the gray cell is disconnected from the rest of the figure. Nevertheless it suffices that each polyomino be horizontally connected.
 Consider first the family of directed plateau polycubes as defined in \cite{CohenSolal,ChamparnaudDJ13} 
 (a polycube is said to be directed if each of its cells can
be reached from a distinguished cell, called root, by a path only made of East, North and  Ahead steps).
The generating function of the number $P_{v,h}$ of directed plateau polycubes of height $h$ and volume $v$ has been computed
in \cite{CohenSolal,ChamparnaudDJ13}:
\begin{equation}\label{plateau}
\displaystyle\sum_{v,h} P_{v,h} p^v t^h = \displaystyle \frac{t \tau(p)}{1 - t p \tau'(p)}
\end{equation}
where $\tau(x) = \sum_{k \geq 1} \frac{x^k}{1-x^k}$ denotes the generating function of the number $\tau(n)$ of divisors of an integer $n$.
%

\begin{wrapfigure}{r}{0.4\textwidth}
 \centering
		\begin{tikzpicture}[scale=0.5]
				
	
\cube{0}{0}{0}\cube{1}{0}{0}
\cube{0}{0}{1}\cube{1}{0}{1}
\cube{0}{0}{2}\cube{1}{0}{2}
\cube{1}{1}{1}\cube{2}{1}{1}\cube{3}{1}{1}
\cube{1}{1}{2}\cube{2}{1}{2}\cube{3}{1}{2}
\cube{2}{2}{1}\cube{3}{2}{1}
			\end{tikzpicture}

		\begin{tikzpicture}[scale=0.5]
				
	\lcube{0}{0}{0}\lcube{0}{0}{1}\lcube{0}{0}{2}
\lcube{0}{1}{1}\lcube{0}{1}{2}
\lcube{0}{2}{1}

	\bcube{0}{0}{0}\bcube{1}{0}{0}\bcube{1}{1}{0}\bcube{2}{1}{0}
	\bcube{3}{1}{0}
\bcube{2}{2}{0}\bcube{3}{2}{0}
				
\tcube{0}{0}{0}\tcube{1}{0}{0}
\tcube{0}{0}{1}\tcube{1}{0}{1}
\tcube{0}{0}{2}\tcube{1}{0}{2}
\tcube{1}{1}{1}\tcube{2}{1}{1}\tcube{3}{1}{1}
\tcube{1}{1}{2}\tcube{2}{1}{2}\tcube{3}{1}{2}
\tcube{2}{2}{1}\tcube{3}{2}{1}
			\end{tikzpicture}
    \caption{A plateau polycube from two horizontally convex directed polyominoes}
	\label{plapol}
	\end{wrapfigure}
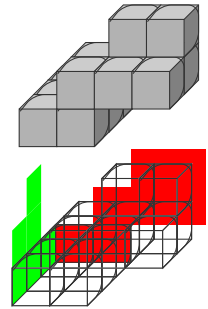

Let us show how to recover (\ref{plateau}) with our method. First we remark that each directed plateau polycube
is obtained from two horizontally convex (\emph{i.e.} each horizontal line meets the polyomino in a single line segment)
 directed (\emph{i.e.} each cell can be reached from $(0,0)$ by movements   up or right one cell, without 
 leaving the polyomino) polyominoes. 
    The generating function of  horizontally convex directed polyominoes of multivolume $[v_1\dots v_h]$ is 
 \[
 \sum_{v_1,\dots,v_h\geq 1}v_1\cdots v_{h-1}x_1^{v_1}\cdots x_h^{v_h}={x_1\cdots x_h\over (1-x_1)^2\cdots(1-x_{h-1})^2(1-x_h)}
 \]
and the convolution yields
\[
  \sum_{v_1 , \dots, v_h \geq 1} v_1 \dots v_{h-1}
   \frac{x_1^{v_1}\dots x_{h}^{v_h}}{(1-x^{v_1})^2\cdots (1-x_{h-1}^{v_{h-1}})^2 (1-x_h^{v_h})}.
\]
Setting $x_1=\dots=x_{h-1}=p$ and $x_h=p\omega$, we obtain the generating function with respect to the
volume ($p$) and the volume of the highest stratum ($\omega$):
$
 \sum_{v_1 , \dots, v_h \geq 1} v_1 \dots v_{h-1}
   \frac{p^{v_1+\cdots+v_{h-1}+v_h}\omega^{v_h}}{(1-p^{v_1})^2\cdots (1-p^{v_{h-1}})^2 (1-(p\omega)^{v_h})} = 
   \left( p\tau'(p) \right)^{h-1} \tau(p\omega).
$
We recover  (\ref{plateau}) by summing over $h$ and setting $\omega=1$ in:
\begin{equation}\label{platomeg}
\sum_{h\geq1}t^h\left( p\tau'(p) \right)^{h-1} \tau(p\omega)={t\tau(p\omega)\over 1-tp\tau'(p)}.
\end{equation}

By a similar reasoning, (\ref{plateau}) can be generalized in dimension $d+1$. If we denote by $P_{v,h}^{(d+1)}$ the number of directed plateau $(d+1)$-polycubes of height $h$ and volume $v$, we can set:
$$\displaystyle{
\sum_{v,h} P_{v,h}^{(d+1)} p^v t^h = \frac{t \tau^{(d+1)}(p)}{1 - t p {\tau^{(d+1)}}'(p)}
}$$
 where   $\tau^{(d+1)}(x) = \left(\frac{x}{1-x}\right)^{*d}$.

The Dirichlet generating function of directed plateau polycubes of height $h$ is $\Xi(s_1-1,\dots,s_{h-1}-1,s_h)^2$ with $\Xi(s_1,\dots,s_h)=\sum_{v_1,\dots,v_h\geq 1}{1\over v_1^{s_1}\cdots v_h^{s_h}}$. 

We give briefly a second example by computing the generating function of the number of all plateau polycubes 
(without constraint). Plateau polycubes of height $h$ are in  one-to-one correspondence 
with pairs of horizontally convex
polyominoes of height $h$. The generating function of   horizontally convex
polyominoes of height $h$  is
\[
\displaystyle\sum_{v_1,\dots,v_h\geq1}(v_1+v_2-1)\cdots(v_{h-1}+v_h-1)x_1^{v_1}\cdots x_h^{v_h}
=\alpha_h(x_1,{x_2\over x_1},\dots,{x_hx_{h-2}\cdots\over x_{h-1}x_{h-3}\cdots})\]
with $\alpha_h(t_1,\dots,t_h)=t_1^2\dots t_{h-1}^2{d\over dt_1}\dots {d\over dt_{h-1}}{t_1\cdots t_h\over
(1-t_1)(1-t_1t_2)(1-t_2t_3)\cdots(1-t_{h-1}t_h)}$.
  Using the same method as above, 
  we obtain the following formula for the generating function of the number of plateau polycubes
   counted by height and volume  :
\[\begin{array}{l}
\displaystyle\sum_{h\geq 1}t^h\sum_{n_1,\dots,n_h}(n_1+n_2-1)\cdots(n_{h-1}+n_h-1)\alpha_h(p^{n_1},p^{n_2-n_1},\dots,p^{n_h-n_{h-1}+n_{h-2}-\cdots})=\\
\displaystyle t \left( \sum_{n_1,m_1\geq 1} p^{n_1 m_1} \left( 1 + t \sum_{n_2,m_2 \geq 1} (n_1+n_2-1)(m_1+m_2-1) p^{n_2m_2} 
 \left( 1 + t \sum_{} \dots \right) \right) \right).
\end{array}
\]

%% file: concl.tex
\section{Conclusion}
In this paper we have investigated the enumeration of pyramids and espaliers
in connection with the Lambert transform. In the last section, we explain that our method
could be used to count other families of polycubes (or more general $3$-dimensional objects).
In particular, we gave  two expressions for the generating function of the plateau polycubes.
Although these expressions are not closed, they can be used to enumerate this family far enough.
Nevertheless, the underlying combinatorics remains to be understood. For instance, a straightforward examination
of the generating function of the horizontally convex polyominoes (see \cite[p. 153]{Wilf}), 
$
\sum_h t^h\alpha_h(p,1,p,1\dots) ={pt(1-p)^3\over (1-p)^4-pt(1-p-p^2+p^3+p^2t)} 
$, 
reveals interesting connections with Delannoy numbers $D_{n-k,k}$ which count 
  the number of lattice paths from $(0,0)$ to $(n,k)$ using steps $(1,0)$, $(0,1)$, $(1,1)$
  (see sequence \href{http://oeis.org/A008288}{A008288} of \cite{Sloane}) and whose generating function is
  $\sum_{n,k}D_{n,k}x^ny^k=(1-x-y-xy)^{-1}$. More precisely, $\alpha_h(p,1,p,1,\dots)={p^h\over (1-p)^{2h-1}}
  \sum_{k=0}^{h-1}D_{h-k-1,k}p^k$. 
  It is natural to ask the question of analogous connections for higher dimensions.\\
One of the tricks allowing to enumerate polyominoes is to consider the statistic of the area of the highest stratum. 
For instance, the generating function of parallelogram polyominoes is deduced from functional equations involving the variable
associated to this statistic (see \emph{e.g.} \cite[Example IX.14 p. 660]{Fla}).
Formula (\ref{platomeg}) shows that this strategy is compatible with our method at least in certain cases. 
It should be interesting to see if one can adapt the ``adding a slice'' method for computing
 functional equations to the generating functions of some families of polycubes obtained from two  polyominoes. Perhaps this method could be adapted by introducing new variables for the width and the length of the highest plateau as suggested by Equation (\ref{platomeg}).

%% file: ref.tex

%% file: polycubes.bbl
\begin{thebibliography}{}

\end{thebibliography}


\begin{thebibliography}{8}

\bibitem{cresson}{ J. Cresson, S. Fischler T. Rivoal}, S\'eries hyperg\'eom\'etriques multiples et polyz\^etas, \url{http://hal.archives-ouvertes.fr/hal-00101376/PDF/CFRalgo-HAL.pdf}.

\bibitem{oai:arXiv.org:math-ph/0212015}{ M. Baake and U. Grimm}, Combinatorial problems of (quasi-)crystallography, in: Quasicrystals - Structure and Physical Properties, 160--171, ed: Wiley-VCH (2002).

\bibitem{Costermans:2009:NAM}{ Christian Costermans and Hoang Ngoc Minh}, Noncommutative algebra, multiple harmonic sums and applications in discrete probability, Journal of Symbolic Computation, vol 44, nb 7, 801--817, 2009.

 
\bibitem{AB06}{ G. Aleksandrowicz, G. Barequet}, Counting $d$-dimensional polycubes and nonrectangular planar polyominoes, in Proc: 12th Ann. Int. Computing Combinatorics Conf., Taipei, Taiwan in: Lecture Notes in Computer Science, vol. 4112, Spinger-Verlag (2006), 418-427.

\bibitem{Jea10}{ H. Jeanne}, Langages g{\'e}om{\'e}triques et polycubes, PhD Thesis, Universit\'e de Rouen, France, 2010.

\bibitem{Bres99}{ Bressoud, D.M.}, Proofs and Confirmations: The Story of the Alternating Sign Matrix Conjecture, in: Spectrum Series, Cambridge University Press (1999).

\bibitem{Stan01}{ Stanley, R.P.}, Enumerative combinatorics. 2, Chapter 7, Cambridge Studies in Advanced Math, Wadsworth \& Brooks/Cole Advanced Books \& Software (2001).

\bibitem{CohenSolal}{ J.-M. Champarnaud, Q. Cohen-Solal, J.-P. Dubernard and H. Jeanne}, Enumeration of Specific Classes of Polycubes (2013).

\bibitem{BK72}{ E.A. Bender, D.E. Knuth}, Enumeration of Plane Partitions, J. Combin. Theory Ser. A. 13 (1972), 40-54.

\bibitem{Bo96}{ M. Bousquet-M\'elou}, A method for the enumeration of various classes of column-convex polygons, Discrete Math., vol 154 (1996), 4-25.


\bibitem{BH57}{ S.R. Broadbent, J.M. Hammersley}, Percolation processes; I. 
Crystals and mazes, Proc. Camb. Philos. Soc. 53 (1957),629-641.

\bibitem{CDJ09}{ J.-M. Champarnaud, J.-Ph. Dubernard, H. Jeanne}, An efficient algorithm to test whether a binary and prolongeable language is geometrical, 
Int. J. Found. Comput. Sci.Vol. 20(4) (2009),763--774.


\bibitem{CLP98}{H. Cohn, M. Larsen, J. Propp}, The Shape of a Typical Boxed Plane Partition, New York J. Math. 4 (1998), 137-166.

\bibitem{Fer04}{S. Feretic}, A $q$-enumeration of convex polyominoes by festoon approach, Theoret. Comput. Sci., 319(1-3) (2004),333-356.

\bibitem{Fla}{Ph. Flajolet, R. Sedgewick}, Analytic Combinatorics, Cambridge University Press (2009). 

\bibitem{LG08}{G. Largeteau, D. Geniet}, Quantification du taux d'invalidit\'e d'applications temps-r\'eel \`a contraintes strictes, Tech. et Sci. Inf. 27(5) (2008), 589-625

\bibitem{Lu71}{ W.F. Lunnon}, Counting polyominoes, in: A.O.L. Atkin B.J. 
Birch (Eds), Computers in Number Theory, Academic Press, London (1971),347-372. 

\bibitem{ChamparnaudDJ13}{ J.-M. Champarnaud, J.-P. Dubernard and H. Jeanne}, A generic method for the enumeration of various classes of directed polycubes, Discrete Mathematics {\&} Theoretical Computer Science, vol: 15, nb: 1 (183-200) (2013).


\bibitem{Sloane}{OEIS Foundation Inc.}(2011), The On-Line Encyclopedia of Integer Sequences, http://oeis.org.
\bibitem{Temp56}{ H. N. V. Temperley}, Combinatorial problems suggested by the statistical machanics of domains and rubber-like molecules,Phys. Rev. 103 (1956), 1-16. 
\bibitem{Wilf} {H. S. Wilf}, generatingfunctionology, Acad press Inc (1994).


\end{thebibliography}
